\documentclass[a4paper,11pt]{amsart}
\usepackage[english]{babel}
\usepackage[latin1]{inputenc}
\usepackage{amssymb}
\usepackage{amscd}
\usepackage{mathtools}
\usepackage[paper=a4paper,left=3cm,right=3cm,top=2cm,bottom=2.0cm]{geometry}
\usepackage{hyperref}
\usepackage{tikz}
\usepackage{tikz-cd}

\newtheorem{thm}{Theorem}[section]
\newtheorem*{thm-non}{Theorem}

\theoremstyle{definition}

\theoremstyle{remark}
\newtheorem{rem}[thm]{Remark}

\begin{document}

\title{The cubo-cubic transformation and K3 surfaces}

\author{Fabian Reede}
\address{Institut f\"ur Algebraische Geometrie, Leibniz Universit\"at Hannover, Welfengarten 1, 30167 Hannover, Germany}
\email{reede@math.uni-hannover.de}

\keywords{Cremona transformations, K3 surfaces, determinantal hypersurfaces}

\subjclass[2010]{Primary: 14E05; Secondary: 14E07, 14J28, 14M12}

\begin{abstract}
In this note we observe that the Cremona transformation in Oguiso's example of Cremona isomorphic but not projectively equivalent quartic K3 surfaces in $\mathbb{P}^3$ is the classical cubo-cubic transformation of $\mathbb{P}^3$.
\end{abstract}
\maketitle

\section*{Introduction}
Let $D_1$ and $D_2$ be two smooth hypersurfaces of degree $d$ in $\mathbb{P}^n$ which are isomorphic as abstract varieties. It is then a natural question if an isomorphism between them can be obtained by restricting an automorphism or a Cremona transformation of the ambient $\mathbb{P}^n$.\\

We say $D_1$ and $D_2$ are projectively equivalent if there is an automorphism $g:\mathbb{P}^n \rightarrow \mathbb{P}^n$ which restricts to an isomorphism $D_1\rightarrow D_2$. We say $D_1$ and $D_2$ are Cremona isomorphic if there is a Cremona transformation $f:\mathbb{P}^n \dashrightarrow \mathbb{P}^n$ such that the restriction extends to an isomorphism $D_1\rightarrow D_2$.\\  

Assume $n\geq 3$ and $(n,d)\neq(3,4)$, then by a result of Matsumura and Monsky,  \cite{matsu}, $D_1$ and $D_2$ are projectively equivalent (so especially Cremona isomorphic) if they are isomorphic as abstract varieties.

Oguiso showed in \cite{oguiso} that in the special case $(n,d)=(3,4)$ this result is no longer true, by giving an example of two smooth quartic K3 surfaces in $\mathbb{P}^3$ which are Cremona isomorphic but which are not projectively equivalent.\\

The main observation of this note is that the Cremona transformation $f:\mathbb{P}^3 \dashrightarrow \mathbb{P}^3$ in Oguiso's example is the classical cubo-cubic transformation of $\mathbb{P}^3$, see \cite{noether} or \cite{semp}.\\

In this note we work over the field $\mathbb{C}$ and the note consists of three sections. In the first section we recall two constructions of the cubo-cubic Cremona transformation. The second section contains a summary of Oguiso's results. In the third section we combine the first two sections and show how to understand Oguiso's example in terms of the cubo-cubic transformation.

\section{The cubo-cubic transformation}

Let $C$ be a general smooth irreducible curve of genus 3 and degree 6 in $\mathbb{P}^3$. The homogeneous ideal of $C$ is generated by four cubic polynomials $f_i\in \mathbb{C}[\textit{\textbf{x}}]$, which are the $3\times 3$ minors of a general $3\times 4$ matrix $A(\textit{\textbf{x}})$ of linear forms, see e.g. \cite[Exemples 2]{elling}, here $\textit{\textbf{x}}=\left[ x_1,x_2,x_3,x_4\right]$ are the homogeneous coordinates on $\mathbb{P}^3$. These polynomials define a rational map:
\begin{equation*}
\begin{tikzcd}
\varphi: \mathbb{P}^3 \arrow[dashrightarrow]{r} & \mathbb{P}^3. 
\end{tikzcd}
\end{equation*}
It follows from \cite[Theorem 7.2.4]{dolga} that this map is in fact birational with base scheme $C$ and has multidegree $(3,3)$. Furthermore the map it is resolved by the graph $\Gamma_{\varphi}$ of $\varphi$ in $\mathbb{P}^3\times\mathbb{P}^3$, that is we have the following diagram:
\begin{equation*}
\begin{tikzcd}[column sep=small]
& \Gamma_{\varphi} \arrow[rightarrow]{dl}[swap]{\pi_1} \arrow[rightarrow]{dr}{\pi_2} & \\
\mathbb{P}^3 \arrow[dashrightarrow]{rr}{\varphi} &&  \mathbb{P}^3
\end{tikzcd}
\end{equation*} 
where $\pi_i: \Gamma_{\varphi} \rightarrow \mathbb{P}^3$ is the $i$-th projection and $\varphi=\pi_2\circ \pi_1^{-1} $.

The graph $\Gamma_{\varphi}$ is the intersection of three general divisors of bidegree $(1,1)$ in $\mathbb{P}^3\times\mathbb{P}^3$. These divisors are defined by the 3 rows of the matrix $A(\textit{\textbf{x}})$. As the multidegree is $(3,3)$ the birational map $\varphi$ is called cubo-cubic transformation.\\

An equivalent definition of the cubo-cubic transformation is the following, see \cite{katz} or \cite[Example 7.2.6]{dolga}: let $C$ be a general smooth irreducible curve of genus 3 and degree 6 in $\mathbb{P}^3$ and denote the blow up of $C$ in $\mathbb{P}^3$ by $X$. Then $\text{Pic}(X)=\mathbb{Z}\left[H\right]\oplus\mathbb{Z}\left[E\right]$, where $H$ is the pullback of a hyperplane in $\mathbb{P}^3$ to $X$ via the blow up morphism $\sigma: X \rightarrow \mathbb{P}^3$ and $E$ is the exceptional divisor of $\sigma$. The linear system $|3H-E|$ defines a morphism 
\begin{equation*}
\Psi: X\rightarrow \mathbb{P}^3
\end{equation*}
which can be shown to be of degree 1, hence $\Psi$ is birational. The morphism $\Psi$ contracts the strict transform $F$ of the trisecant surface $Sec_3(C)$ with respect to $\sigma$ to a curve $C'\subset \mathbb{P}^3$ with $C\cong C'$. 

We get a birational $\phi$ map from the following diagram:
\begin{equation*}
\begin{tikzcd}[column sep=small]
& X \arrow[rightarrow]{dl}[swap]{\sigma} \arrow[rightarrow]{dr}{\Psi} & \\
\mathbb{P}^3 \arrow[dashrightarrow]{rr}{\phi} &&  \mathbb{P}^3
\end{tikzcd}
\end{equation*}
We can see that $\phi=\varphi$ by noting that the graph $\Gamma_{\varphi}$, i.e. the intersection of the three divisors of bidegree $(1,1)$ in $\mathbb{P}^3\times\mathbb{P}^3$, is the same as the blow-up of the curve $C$, see for example \cite[H.29, Lemma E.1]{coates} or \cite[7.1.3]{dolga}.
\begin{rem}
The cubo-cubic Cremona transformation and most of its properties were already known to Max Noether, see \cite[§2]{noether}.
\end{rem}
\begin{rem}
The cubo-cubic transformation is special in the sense that it is the only non-trivial Cremona transformation of $\mathbb{P}^3$ that is resolved by just one blow up along a smooth curve, see \cite[Proposition 2.1]{katz}.
\end{rem}
\section{Oguiso's example}
In \cite[Theorem 1.5.]{oguiso} Oguiso constructs two quartic K3 surfaces $S_1$ and $S_2$ in $\mathbb{P}^3$ and a Cremona transformation $\tau: \mathbb{P}^3\dashrightarrow \mathbb{P}^3$ such that $\tau$ restricts to a birational map $S_1 \dashrightarrow S_2$. This map must be an isomorphism as the canonical divisor of a K3 surface is nef. But these two surfaces are not projectively equivalent in $\mathbb{P}^3$, that is there is no $g\in \text{Aut}(\mathbb{P}^3)=\text{PGL}(4,\mathbb{C})$ with $g(S_1)=S_2$.

The K3 surfaces are constructed as follows: pick three general divisors $Q_1$, $Q_2$ and $Q_3$ of bidegree $(1,1)$ in $\mathbb{P}^3\times\mathbb{P}^3$, then 
\begin{equation*}
V:=Q_1\cap Q_2\cap Q_3
\end{equation*}
is a smooth Fano threefold, birational to $\mathbb{P}^3$ in two different ways, given by restricting the projections $\pi_i: \mathbb{P}^3\times\mathbb{P}^3 \rightarrow \mathbb{P}^3$ to $V$, $i=1,2$. Call these morphisms $p_1$ and $p_2$.  Oguiso defines the birational map $\tau:\mathbb{P}^3 \dashrightarrow\mathbb{P}^3$ via:
\begin{equation*}
\tau:=p_2\circ p_1^{-1}. 
\end{equation*} 
This defines a non-trivial Cremona transformation of $\mathbb{P}^3$, \cite[Theorem 1.5.(2)]{oguiso}.

For a very general divisor $Q_4$ of bidegree $(1,1)$, avoiding the Noether-Lefschetz locus, define:
\begin{equation*}
S:=Q_1\cap Q_2\cap Q_3\cap Q_4.
\end{equation*}
Then $S$ is a smooth K3 surface with $\text{NS}(S)=\mathbb{Z}h_1\oplus\mathbb{Z}h_2$ with intersection matrix
\begin{equation*}
((h_i,h_j))_{i,j}=\begin{pmatrix} 4 & 6 \\ 6 & 4\end{pmatrix}.
\end{equation*}
The restrictions of the $p_i$ for $i=1,2$, i.e. ${p_i}_{|S}:S\rightarrow p_i(S)$, turn out to be isomorphisms given by the complete linear systems $|h_i|$. We get two quartic K3 surfaces $S_i:=p_i(S)\subset \mathbb{P}^3$.

The main result of Oguiso regarding these two surfaces is:

\begin{thm}\cite[Theorem 1.5.(3)]{oguiso}\label{main}
The K3 surfaces $S_1$ and $S_2$ are Cremona isomorphic but not projectively equivalent in $\mathbb{P}^3$.
\end{thm}

\begin{rem}
The reason why $S_1$ and $S_2$ are not projectively equivalent is that there is no $f\in \text{Aut}(S)$ such that $f^{*}h_1=h_2$, see \cite[Proposition 6.2., Lemma 6.4.]{oguiso}. 
\end{rem}

The whole construction is captured in the following diagram:
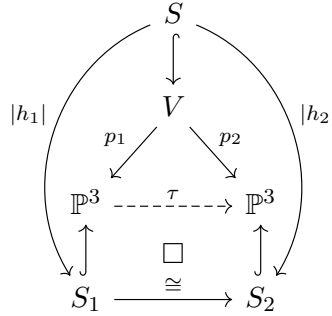
\begin{equation*}
\begin{tikzcd}[column sep=small]
& S \arrow[swap, dddl,bend right=50,"|h_1|"]\arrow[dddr,bend left=50,"|h_2|"]\arrow[hookrightarrow]{d}\\
& V \arrow[rightarrow]{dl}[swap]{p_1} \arrow[rightarrow]{dr}{p_2} & \\
\mathbb{P}^3 \arrow[dashrightarrow]{rr}{\tau} &&  \mathbb{P}^3  \\
S_1 \arrow[hookrightarrow]{u} \arrow{rr}{\cong}\arrow[phantom, "\square"]{urr} && S_2\arrow[hookrightarrow]{u} 
\end{tikzcd}
\end{equation*}

\begin{rem}
The K3 surfaces $S_1$ and $S_2$ are determinantal quartic K3 surfaces and were already known to Cayley, see \cite{cay} and \cite{fest}.
\end{rem}

\section{Cubo-cubic Cremona transformations and K3 surfaces}
We want to bring the previous two sections together and understand Oguiso's example in terms of a cubo-cubic transformation. To do this, pick a general smooth irreducible curve $C$ of genus 3 and degree 6 in $\mathbb{P}^3$, with $3\times 4$ matrix $A(\textit{\textbf{x}})$ and write
\begin{equation*}
A(\textit{\textbf{x}})=x_1A_1+x_2A_2+x_3A_3+x_4A_4
\end{equation*}
with $A_k=(a_{ij}^k)_{i,j}\in \text{Mat}(3,4,\mathbb{C})$. This defines three $4\times 4$ matrices $B_i=(a_{ij}^k)_{j,k}$ for $i=1,2,3$ where $k$ is the index for the columns and $j$ is the index for the rows.

The three divisors $Q_1$, $Q_2$ and $Q_3$ of bidegree $(1,1)$ induced by $A(\textit{\textbf{x}})$ (the equations of the graph of the cubo-cubic transformation induced by $C$) are :
\begin{equation*}
Q_i=\left(\sum\limits_{k,j=1}^4a_{ij}^kx_ky_j=0 \right)\subset \mathbb{P}^3\times\mathbb{P}^3\,\,\,\,\,\,\,(i=1,2,3).
\end{equation*} 
We choose a very general fourth divisor $Q_4$ of bidegree $(1,1)$ as described above:
\begin{equation*}
Q_4=\left(\sum\limits_{k,j=1}^4a_{4j}^kx_ky_j=0 \right)\subset \mathbb{P}^3\times\mathbb{P}^3.
\end{equation*} 
Thus $S=Q_1\cap Q_2\cap Q_3\cap Q_4 \subset \mathbb{P}^3\times\mathbb{P}^3$ is given by
\begin{equation*}
S=\left\lbrace (\textbf{x},\textbf{y})\in\mathbb{P}^3\times\mathbb{P}^3\,|\, M(\textit{\textbf{x}})\cdot \textit{\textbf{y}}^t=0^t \right\rbrace  
\end{equation*}
for the $4\times 4$ matrix $M(\textit{\textbf{x}})=\left( m_{ij}(\textit{\textbf{x}})\right) _{i,j}$ with entries:
\begin{equation*}
m_{ij}(\textit{\textbf{x}})=a_{ij}^1x_1+a_{ij}^2x_2+a_{ij}^3x_3+a_{ij}^4x_4.
\end{equation*}
The observation
\begin{equation}\label{obs}
\sum\limits_{k,j=1}^4a_{ij}^kx_ky_j=\sum\limits_{j=1}^4\left(\sum\limits_{k=1}^4 a_{ij}^kx_k \right)y_j=\sum\limits_{k=1}^4\left(\sum\limits_{j=1}^4 a_{ij}^ky_j \right)x_k,
\end{equation}
implies the following identity: 
\begin{equation*}
M(\textit{\textbf{x}})\cdot \textit{\textbf{y}}^t=N(\textit{\textbf{y}})\cdot \textit{\textbf{x}}^t
\end{equation*}
with the $4\times 4$ matrix $N(\textit{\textbf{y}})=\left(n_{ik}(\textit{\textbf{y}}) \right)_{i,k}$ given by
\begin{equation*}
n_{ik}(\textit{\textbf{y}})=a_{i1}^ky_1+a_{i2}^ky_2+a_{i3}^ky_3+a_{i4}^ky_4.
\end{equation*}
Oguiso proved in \cite[Proposition 4.1.]{oguiso} that 
\begin{equation*}
S_1=\left\lbrace\, \textbf{x}\in \mathbb{P}^3\,|\, \text{det}(M(\textit{\textbf{x}}))=0\,\right\rbrace \,\,\text{and}\,\, S_2=\left\lbrace \textbf{y}\in \mathbb{P}^3\,|\, \,\text{det}(N(\textit{\textbf{y}}))=0\,\right\rbrace .
\end{equation*}
By construction the first three rows of $M(\textit{\textbf{x}})$ are those of $A(\textit{\textbf{x}})$. Using the Laplace expansion with respect to the last row shows that for every $\textit{\textbf{x}}\in C$ we have $\text{det}(M(\textit{\textbf{x}}))=0$ and hence $C\subset S_1$.

Similarly the first three rows of $N(\textit{\textbf{y}})$ define a general $3\times 4$ matrix $A'(\textit{\textbf{y}})$ whose $3\times 3$ minors give rise to a genus 3 and degree 6 curve $C'$ with $C'\subset S_2$.

The blow up $X$ of $C$ in $\mathbb{P}^3$ is given by $Q_1\cap Q_2\cap Q_3$, which by (\ref{obs}) is also the blow up of $C'$ in $\mathbb{P}^3$ and $C\cong C'$ by \cite[Proposition 1.3.]{katz}.

Denote by $\widetilde{S}_1$ the strict transform of $S_1$ in the blow up $X$ of $C$ in $\mathbb{P}^3$ and by $\widetilde{S}_2$ the strict transform of $S_2$ in the blow up of $C'$.
Using Oguiso's results and the fact that $C$ and $C'$ are smooth we get 
\begin{equation*}
S_1\cong \widetilde{S}_1=S=\widetilde{S}_2\cong S_2.
\end{equation*}

Finally we look at the cubo-cubic Cremona transformation given by the curve $C$:
\begin{equation*}
\begin{tikzcd}[column sep=small]
& X \arrow[rightarrow]{dl}[swap]{\sigma} \arrow[rightarrow]{dr}{\Psi} & \\
\mathbb{P}^3 \arrow[dashrightarrow]{rr}{\phi} &&  \mathbb{P}^3
\end{tikzcd}
\end{equation*}  

By the previous results the strict transform $\widetilde{S}_1$ with respect to $\sigma$ of the K3 surface $S_1$ in $X$ is the K3 surface $S$ and it is also equal to the strict transform $\widetilde{S}_2$ with respect to $\Psi$ of the K3 surface $S_2$.

This implies that the birational map $\phi$ restricts to a birational map between $S_1$ and $S_2$ which extends to an isomorphism, because $S_1$ and $S_2$ are K3 surfaces, which have nef canonical divisors. Hence $S_1$ and $S_2$ are Cremona isomorphic in $\mathbb{P}^3$ and the Cremona isomorphism is induced by the cubo-cubic transformation given by the curve $C$. But by Theorem \ref{main}, they are not projectively equivalent in $\mathbb{P}^3$.

\begin{rem}
The construction of smooth determinantal quartic surfaces in $\mathbb{P}^3$ containing a curve of genus $3$ and degree $6$ can also be found in Beauville's paper \cite[6.7]{beau}. In fact Beauville proves: a smooth quartic surface in $\mathbb{P}^3$ is determinantal if and only if contains a nonhyperelliptic curve of genus 3 embedded in $\mathbb{P}^3$ by a linear system of degree 6, see \cite[Corollary 6.6]{beau}.
\end{rem}

\bibliography{Artikel}

\begin{thebibliography}{FGvGvL13}

\bibitem[Bea00]{beau}
Arnaud Beauville.
\newblock Determinantal hypersurfaces.
\newblock {\em Michigan Math. J.}, 48:39--64, 2000.

\bibitem[Cay71]{cay}
Arthur Cayley.
\newblock A {M}emoir on {Q}uartic {S}urfaces.
\newblock {\em Proc. Lond. Math. Soc.}, 3:19--69, 1869/71.

\bibitem[CCGK16]{coates}
Tom Coates, Alessio Corti, Sergey Galkin, and Alexander Kasprzyk.
\newblock Quantum periods for 3-dimensional {F}ano manifolds.
\newblock {\em Geom. Topol.}, 20(1):103--256, 2016.

\bibitem[Dol12]{dolga}
Igor~V. Dolgachev.
\newblock {\em Classical algebraic geometry - A modern view}.
\newblock Cambridge University Press, Cambridge, 2012.

\bibitem[Ell75]{elling}
Geir Ellingsrud.
\newblock Sur le sch\'{e}ma de {H}ilbert des vari\'{e}t\'{e}s de codimension
  {$2$} dans {${\bf P}\sp{e}$} \`a c\^{o}ne de {C}ohen-{M}acaulay.
\newblock {\em Ann. Sci. \'{E}cole Norm. Sup. (4)}, 8(4):423--431, 1975.

\bibitem[FGvGvL13]{fest}
Dino Festi, Alice Garbagnati, Bert van Geemen, and Ronald van Luijk.
\newblock The {C}ayley-{O}guiso automorphism of positive entropy on a {K}3
  surface.
\newblock {\em J. Mod. Dyn.}, 7(1):75--97, 2013.

\bibitem[Kat87]{katz}
Sheldon Katz.
\newblock The cubo-cubic transformation of {${\bf P}^3$} is very special.
\newblock {\em Math. Z.}, 195(2):255--257, 1987.

\bibitem[MM64]{matsu}
Hideyuki Matsumura and Paul Monsky.
\newblock On the automorphisms of hypersurfaces.
\newblock {\em J. Math. Kyoto Univ.}, 3:347--361, 1963/1964.

\bibitem[Noe71]{noether}
Max Noether.
\newblock Ueber die eindeutigen {R}aumtransformationen, insbesondere in ihrer
  {A}nwendung auf die {A}bbildung algebraischer {F}l\"{a}chen.
\newblock {\em Math. Ann.}, 3(4):547--580, 1871.

\bibitem[Ogu17]{oguiso}
Keiji Oguiso.
\newblock Isomorphic quartic {K}3 surfaces in the view of {C}remona and
  projective transformations.
\newblock {\em Taiwanese J. Math.}, 21(3):671--688, 2017.

\bibitem[SR49]{semp}
John~G. Semple and Leonard Roth.
\newblock {\em Introduction to {A}lgebraic {G}eometry}.
\newblock Oxford, at the Clarendon Press, 1949.

\end{thebibliography}
\bibliographystyle{alpha}

\end{document}